\catcode`\^^Z=9
\catcode`\^^M=10
\def\RingvGN{1.1}

\def\calGgen{2.1}
\def\ActNode{3.1}

\def\IsoStabs{3.4}

\def\DefCrep{4.1}
\def\LocCrep{4.2}
\def\Ito{4.3}
\def\Diffinv{5.1}

\def\Fort{5.3}

\def\picGN{6.1}

\def\DecIrr{6.4}
\def\somPp{7.1}

\output={\if N\header\headline={\hfill}\fi
\plainoutput\global\let\header=Y}
\magnification\magstep1
\tolerance = 500
\hsize=14.4true cm
\vsize=22.5true cm
\parindent=6true mm\overfullrule=2pt
\newcount\kapnum \kapnum=0
\newcount\parnum \parnum=0
\newcount\procnum \procnum=0
\newcount\nicknum \nicknum=1
\font\ninett=cmtt9

\font\ninebf=cmbx9

\font\sixbf=cmbx6
\font\ninesl=cmsl9

\font\nineit=cmti9

\font\ninerm=cmr9

\font\sixrm=cmr6
\font\ninei=cmmi9
\font\eighti=cmmi8
\font\sixi=cmmi6
\skewchar\ninei='177 \skewchar\eighti='177 \skewchar\sixi='177
\font\ninesy=cmsy9
\font\eightsy=cmsy8
\font\sixsy=cmsy6
\skewchar\ninesy='60 \skewchar\eightsy='60 \skewchar\sixsy='60
\font\titelfont=cmr10 scaled 1440
\font\paragratit=cmbx10 scaled 1200

\font\name=cmcsc10
\font\emph=cmbxti10

\font\tenmsbm=msbm10
%\font\ninemsbm=msbm9
\font\sevenmsbm=msbm7
%\font\sixmsbm=msbm6
%\font\fivemsbm=msbm5
%\textfont\extsym=\tenmsbm
%\scriptfont\extsym=\sevenmsbm
%\scriptscriptfont\extsym=\fivemsbm
%

%

\font\Got=eufm7
\font\teneufm=eufm10
\font\seveneufm=eufm7
\font\fiveeufm=eufm5
\newfam\eufmfam
\textfont\eufmfam=\teneufm
\scriptfont\eufmfam=\seveneufm
\scriptscriptfont\eufmfam=\fiveeufm

\font\tenmsam=msam10
\font\sevenmsam=msam7
\font\fivemsam=msam5
\newfam\msamfam
\textfont\msamfam=\tenmsam
\scriptfont\msamfam=\sevenmsam
\scriptscriptfont\msamfam=\fivemsam
\font\tenmsbm=msbm10
\font\sevenmsbm=msbm7
\font\fivemsbm=msbm5
\newfam\msbmfam
\textfont\msbmfam=\tenmsbm
\scriptfont\msbmfam=\sevenmsbm
\scriptscriptfont\msbmfam=\fivemsbm
\def\Bbb#1{{\fam\msbmfam\relax#1}}
\def\cz{{\kern0.4pt\Bbb C\kern0.7pt}
}
\def\ez{{\kern0.4pt\Bbb E\kern0.7pt}
}
\def\fz{{\kern0.4pt\Bbb F\kern0.3pt}}
\def\gz{{\kern0.4pt\Bbb Z\kern0.7pt}}
\def\hz{{\kern0.4pt\Bbb H\kern0.7pt}
}
\def\kz{{\kern0.4pt\Bbb K\kern0.7pt}
}
\def\nz{{\kern0.4pt\Bbb N\kern0.7pt}
}
\def\oz{{\kern0.4pt\Bbb O\kern0.7pt}
}
\def\rz{{\kern0.4pt\Bbb R\kern0.7pt}
}
\def\sz{{\kern0.4pt\Bbb S\kern0.7pt}
}
\def\pz{{\kern0.4pt\Bbb P\kern0.7pt}
}
\def\qz{{\kern0.4pt\Bbb Q\kern0.7pt}
}
\newskip\ttglue
\def\ninepoint{\def\rm{\fam0\ninerm}%
  \textfont0=\ninerm \scriptfont0=\sixrm \scriptscriptfont0=\fiverm
  \textfont1=\ninei \scriptfont1=\sixi \scriptscriptfont1=\fivei
  \textfont2=\ninesy \scriptfont2=\sixsy \scriptscriptfont2=\fivesy
  \textfont3=\tenex \scriptfont3=\tenex \scriptscriptfont3=\tenex
  \def\it{\fam\itfam\nineit}%
  \textfont\itfam=\nineit
  \def\sl{\fam\slfam\ninesl}%
  \textfont\slfam=\ninesl
  \def\bf{\fam\bffam\ninebf}%
  \textfont\bffam=\ninebf \scriptfont\bffam=\sixbf
   \scriptscriptfont\bffam=\fivebf
  \def\tt{\fam\ttfam\ninett}%
  \textfont\ttfam=\ninett
  \tt \ttglue=.5em plus.25em minus.15em
  \normalbaselineskip=11pt
  \font\name=cmcsc9
  \let\sc=\sevenrm
  \let\big=\ninebig
  \setbox\strutbox=\hbox{\vrule height8pt depth3pt width0pt}%
  \normalbaselines\rm
  \def\sl{\it}}

\headline={\ifodd\pageno\rightheadline\else\leftheadline\fi}
\def\rightheadline{\ninepoint Paragraphen"uberschrift\hfill\folio}
\def\leftheadline{\ninepoint\folio\hfill Chapter"uberschrift}
\let\header=Y
\def\titel#1{\need 9cm \vskip 2truecm
\parnum=0\global\advance \kapnum by 1
{\baselineskip=16pt\lineskip=16pt\rightskip0pt
plus4em\spaceskip.3333em\xspaceskip.5em\pretolerance=10000\noindent
\titelfont Chapter \uppercase\expandafter{\romannumeral\kapnum}.
#1\vskip2true cm}\def\leftheadline{\ninepoint
\folio\hfill Chapter \uppercase\expandafter{\romannumeral\kapnum}.
#1}\let\header=N
}
\def\Titel#1{\need 9cm \vskip 2truecm
\global\advance \kapnum by 1
{\baselineskip=16pt\lineskip=16pt\rightskip0pt
plus4em\spaceskip.3333em\xspaceskip.5em\pretolerance=10000\noindent
\titelfont\uppercase\expandafter{\romannumeral\kapnum}.
#1\vskip2true cm}\def\leftheadline{\ninepoint
\folio\hfill\uppercase\expandafter{\romannumeral\kapnum}.
#1}\let\header=N
}
\def\need#1cm {\par\dimen0=\pagetotal\ifdim\dimen0<\vsize
\global\advance\dimen0by#1 true cm
\ifdim\dimen0>\vsize\vfil\eject\noindent\fi\fi}
\def\neupara#1{\par\penalty-2000
\procnum=0\global\advance\parnum by 1
\vskip1cm\noindent{\paragratit \the\parnum. #1}%
\def\rightheadline{\ninepoint\S\the\parnum.\ #1\hfill \folio}%
\vskip 8mm\noindent}
\def\Proclaim #1 #2\finishproclaim {\bigbreak\noindent
{\bf#1\unskip{}. }{\it#2}\medbreak\noindent}
%
%\llap{#3\ }
\gdef\proclaim #1 #2 #3\finishproclaim {\bigbreak\noindent%
\global\advance\procnum by 1
{%
{\relax\ifodd \nicknum
\hbox to 0pt{\vrule depth 0pt height0pt width\hsize
   \quad \ninett#3\hss}\else {}\fi}%
\bf\the\parnum.\the\procnum\ #1\unskip{}. }
{\it#2}
\immediate\write\num{\string\def
 \expandafter\string\csname#3\endcsname
 {\the\parnum.\the\procnum}}
\medbreak\noindent}
\newcount\stunde \newcount\minute \newcount\hilfsvar
\def\uhrzeit{
    \stunde=\the\time \divide \stunde by 60
    \minute=\the\time
    \hilfsvar=\stunde \multiply \hilfsvar by 60
    \advance \minute by -\hilfsvar
    \ifnum\the\stunde<10
    \ifnum\the\minute<10
    0\the\stunde:0\the\minute~Uhr
    \else
    0\the\stunde:\the\minute~Uhr
    \fi
    \else
    \ifnum\the\minute<10
    \the\stunde:0\the\minute~Uhr
    \else
    \the\stunde:\the\minute~Uhr
    \fi
    \fi
    }

\def\calC{{\cal C}} 
 
\def\calG{{\cal G}} \def\calH{{\cal H}}
 
\def\calK{{\cal K}}

 \def\calX{{\cal X}}
 \def\calZ{{\cal Z}}

\def\Gotn{\hbox{\Got n}}

\def\dim{\mathop{\rm dim}\nolimits}

\def\GL{\mathop{\rm GL}\nolimits}

\def\id{\mathop{\rm id}\nolimits}

\def\kernel{\mathop{\rm kernel}\nolimits}

\def\mod{\mathop{\rm mod}\nolimits}

\def\Pic{\mathop{\rm Pic}\nolimits}

\def\proj{\mathop{\rm proj}\nolimits}

\def\SL{\mathop{\rm SL}\nolimits}

\def\Sp{\mathop{\rm Sp}\nolimits}

\def\boxit#1{
  \vbox{\hrule\hbox{\vrule\kern6pt
  \vbox{\kern8pt#1\kern8pt}\kern6pt\vrule}\hrule}}
\def\Boxit#1{
  \vbox{\hrule\hbox{\vrule\kern2pt
  \vbox{\kern2pt#1\kern2pt}\kern2pt\vrule}\hrule}}

\def\smallni{\smallskip\noindent }
\def\medni{\medskip\noindent }
\def\bigni{\bigskip\noindent }
\def\Isom{\mathop{\;{\buildrel \sim\over\longrightarrow }\;}}
\def\lo{\longrightarrow}

\def\loma{\longmapsto}
\def\imag{{\rm i}}
\def\pii{\pi {\rm i}}

\def\set#1{\bigl\{\,#1\,\bigr\}}

\def\square{\hbox{\hbox to 0pt{$\sqcup$\hss}\hbox{$\sqcap$}}}
\def\qed{\ifmmode\square\else{\unskip\nobreak\hfil
\penalty50\hskip3em\null\nobreak\hfil\square
\parfillskip=0pt\finalhyphendemerits=0\endgraf}\fi}
\def\pn{\the\parnum.\the\procnum}
\def\downmapsto{{\buildrel
        {\vbox{\hbox{\hskip.2pt$\scriptstyle-$}}}
        \over{\raise7pt\vbox{\vskip-4pt\hbox{$\textstyle\downarrow$}}}}}
\nopagenumbers
\immediate\newwrite\num
%\nicknum=0  %bewirkt, dass keine Nicknamen am Rand
\font\ttt=cmtt10
\def\hodge#1#2%
{\matrix{&&&1&&&\cr&&0&&0&&\cr&0&&#1&&0&\cr
1&&#2&&#2&&1\cr&0&&#1&&0&\cr&&0&&0&&\cr&&&1&&&\cr}}
\def\t{\vartheta}%
\def\tt#1#2{{\t\Bigl[\matrix{{#1}\cr {#2}}\Bigl]}}%
\font\dunh=cmss10
\let\header=N
\def\PGL{{\rm PGL}}
\def\transpose#1{\kern1pt{^t\kern-3pt#1}}%
\def\pic{{\rm pic}}
\def\Fix{{\rm Fix}}

\immediate\openout\num=calabi6.num
\immediate\newwrite\num\immediate\openout\num=calabi6.num
\def\RAND#1{\vskip0pt\hbox to 0mm{\hss\vtop to 0pt{%
  \raggedright\ninepoint\parindent=0pt%
  \baselineskip=1pt\hsize=2cm #1\vss}}\noindent}
\noindent
\centerline{\titelfont Some Siegel threefolds with a Calabi-Yau model II}%
\def\leftheadline{\ninepoint\folio\hfill
Some Siegel threefolds with a Calabi-Yau model II}%
\def\rightheadline{\ninepoint Introduction\hfill \folio}%
\headline={\ifodd\pageno\rightheadline\else\leftheadline\fi}
\vskip 1.5cm
\leftline{\ttt \hbox to 6cm{Eberhard Freitag\hss}
Riccardo Salvati
Manni  }
  \leftline {\it  \hbox to 6cm{Mathematisches Institut\hss}
Dipartimento di Matematica, }
\leftline {\it  \hbox to 6cm{Im Neuenheimer Feld 288\hss}
Piazzale Aldo Moro, 2}
\leftline {\it  \hbox to 6cm{D69120 Heidelberg\hss}
 I-00185 Roma, Italy. }
\leftline {\ttt \hbox to 6cm{freitag@mathi.uni-heidelberg.de\hss}
salvati@mat.uniroma1.it}
\vskip1cm
\centerline{\paragratit \rm  2010}%
\vskip5mm\noindent%
\let\header=N%
\def\imag{{\rm i}}%
{\paragratit Introduction}%
\medni
In the paper [FSM] we described some Siegel modular threefolds which admit
a Calabi-Yau model. Using a different method we give in this
paper an enlarged list of such varieties. Basic for our method is
the paper [GN] of van Geemen and Nygaard. In this paper they prove
that the complete intersection $\calX$ in $P^7(\cz)$, defined by the equations
$$\eqalign{Y_0^2&=X_0^2+X_1^2+X_2^2+X_3^2,\cr
Y_1^2&=X_0^2-X_1^2+X_2^2-X_3^2,\cr
Y_2^2&=X_0^2+X_1^2-X_2^2-X_3^2,\cr
Y_3^2&=X_0^2-X_1^2-X_2^2+X_3^2\cr}\leqno\qquad\calX:$$
is biholomorphic equivalent to  the Satake compactification of
$\hz_2/\Gamma'$ for
a certain subgroup $\Gamma'\subset\Sp(2,\gz)$. This variety has 96
singularities which correspond to certain zero-dimensional cusps
and these singularities are ordinary double points (nodes).
In the paper [CM] it has
been pointed out that the results of [GN] imply that
a (projective) small resolution of this variety is a
rigid Calabi-Yau manifold $\tilde\calX$.
\smallskip
We describe the basic occurring groups: We use the standard notations,
$M={A\,B\choose C\,D}$:
$$\eqalign{
\Gamma_n[l]&=\kernel (\Sp(n,\gz)\lo\Sp(n,\gz/l\gz)),\cr
\Gamma_n[l,2l]&=\set{M\in\Gamma_n[l];\quad
A{\transpose B}/l\ \hbox{and}\ C{\transpose D}/l\ \hbox{have even diagonal}},\cr
\Gamma_{n,0}[l]&=\{M\in\Gamma_n;\quad C\equiv 0\;\mod\; l\},\cr
\Gamma_{n,0,\vartheta}[l]&=\{M\in\Gamma_n;\quad C\equiv 0\;\mod\; l,
\quad C{\transpose D}/l\ \hbox{has even diagonal}\}.\cr
}$$
The group $\Gamma_{n,0}[l]$ can be extended by the Fricke involution
$$J=J_l=\pmatrix{0&E/\sqrt l\cr-\sqrt lE&0}
\qquad\Bigl(JZ=-{1\over lZ}\Bigr).$$
We denote this extension (of index 2) by
$$\hat\Gamma_{n,0}[l]=\Gamma_{n,0}[l]\cup J\Gamma_{n,0}[l].$$
The group $\Gamma'$, which belongs to van Geemen's and Nygaard's variety
is a subgroup of index two of the group
$$
\Gamma_2[2,4]\cap\Gamma_{2,0,\vartheta}[4]
=\{M\in\Gamma_2[2,4];\quad C\equiv 0\;\mod\; 4,\ C/4\
\hbox{has even diagonal}\},$$
namely
$$\Gamma'=\{M\in \Gamma_2[2,4]\cap\Gamma_{2,0,\vartheta}[4];
\quad \det D\equiv\pm1\;\mod\;8\}.$$
In [FSM] we introduced a certain subgroup $\Gamma_{2,0}[2]_{\Gotn}$
of index two of $\Gamma_{2,0}[2]$ as kernel of a certain
character $\chi_{\Gotn}$. This character is the product
of the unique non-trivial character of the full Siegel modular group
and the character
$$\imag^{\alpha+\beta+ \gamma},\quad C
{\transpose D} = \pmatrix{\alpha&\beta\cr \beta& \gamma},
\qquad\Bigl(M=\pmatrix{A&B\cr C&D}\in\Gamma_{2,0}[2]\Bigr).$$
The group $\Gamma'$ is contained in
$\Gamma_{2,0}[2]_{\Gotn}$ as subgroup
of index 6144.
\smallskip
The character $\chi_{\Gotn}$ extends to a character $\hat\chi_{\Gotn}$
of $\hat\Gamma_{2,0}[2]$ by means of
$$\hat\chi_{\Gotn}(J)=1.$$
We denote its kernel by $\hat\Gamma_{2,0}[2]_{\Gotn}$. This is an
extension of index two of $\Gamma_{2,0}[2]_{\Gotn}$.
The group $\Gamma'$ is contained in
$\hat\Gamma_{2,0}[2]_{\Gotn}$ as subgroup
of index $12\,288$.
The main result of this paper is:
\Proclaim
{Theorem}
{The Siegel modular threefold, which belongs to a group between
$\Gamma'$ and $\hat\Gamma_{2,0}[2]_{\Gotn}$, admits a Calabi-Yau model
in the following weak sense:
There exists a desingularization in the category of complex spaces
of the
Satake compactification which admits a holomorphic three-form without zeros
and whose first Betti number vanishes.}
\finishproclaim
It has been pointed out by van Geemen that it is not always possible to get a
projective model. We want to come back to the question of projectivity in an extra
paper. Here we just mention a result  from [CFS] that can be proved with the
methods of [FSM]:
\Proclaim
{Supplement}
{In the case of a group between $\Gamma_{2,0}[2]_{\Gotn}$ and
$\Gamma_{2,0}[4]\cap\Gamma_2[2]$ one can get a projective
Calabi-Yau manifold.}
\finishproclaim
We will develop a method to compute the Picard number $\pic$
(rank of the group of divisor classes)
and the Euler number for the
intermediate groups. In the projective case they determine the Hodge numbers,
the Picard number  is $\pic=h^{11}$ and the Euler number is
$$e=2(h^{11}-h^{12}).$$
For a small resolution $\tilde\calX$
of the variety of van Geemen and Nygaard $\calX$ they are known
[GN], [CM],
$\pic=32$, $e=64$.  To get it for other groups, one needs the action
of the group $\hat\Gamma_{2,0}[2]_{\Gotn}$  on the Picard group
of the regular locus of $\calX$. We will determine this action
in section 6. The result of this section allow in principle
to compute the numbers for all intermediate groups.
There are thousands of conjugacy classes of intermediate groups.
\smallskip
In section 7 we treat some simple examples, namely all subgroups
of order two of $\hat\Gamma_{2,0}[2]_{\Gotn}/\Gamma'$.
\smallskip
In a subsequent paper we will treat some more  significative cases.
%
%and in section 8 a subgroup of order 3. This group is essentially
%unique since all subgroups of order 3 are conjugated. Actually
%this group of order 3 is contained in a cyclic subgroup of order 6,
%which also leads to an example of a Calabi-Yau manifold.
%\smallskip
%There is one intermediate group of particular interest, namely the group
%$$\Gamma_2[2]\cap\Gamma_{2,0}[4].$$
%This groups contains $\Gamma'$ as subgroup of index 32.
%For this group (and as a consequence for all groups between it
%and $\Gamma_{2,0}[2]_{\Gotn}$) we have a completely different proof which rests
%on the paper [FSM] and gives a very explicit description of the
%Calabi-Yau model, namely:
%\Proclaim
%{Theorem}
%{Let $\tilde X(4)$ be the Igusa desingularization of the Satake compactification
%of $\hz_2/\Gamma[4]$. Then the quotient
%$\tilde X(4)/(\Gamma_2[2]\cap\Gamma_{2,0}[4])$
%admits a desingularization, which is a Calabi-Yau manifold.}
%\finishproclaim
%Since this particular description may be of own interest, we
%add a comment at the end of the paper, how this result follows from [FSM].
%Finally we determine the ring of modular forms for this special group
%and compute the Hodge numbers of a Calabi-Yau model.
\smallskip
We are grateful for helpful disussions
with Bert van Geemen and also for useful comments
to a preliminary version of our paper from Philip Candelas.
He brought our attention to the paper [BH] of Borisov and Hua, where other
examples of complete intersections of 4 quadrics in $P^7(\cz)$ that lead to
Calabi-Yau manifolds with big fundamental groups are described.
\neupara{The variety of van Geemen and Nygaard}%
We recall some basic facts about Siegel modular varieties.
For details we refer to [Fr1].
The symplectic
group $\Sp(n,\rz)$ acts on the Siegel upper half plane
$$\hz_n:=\{Z\in\cz^{(n,n)};\quad
Z=\transpose Z=X+\imag Y,\ Y>0\ \hbox{\rm (positive
definite)}\}$$
by means of the formula
$MZ=(AZ+B)(CZ+D)^{-1}$. For any subgroup $\Gamma\subset\Sp(n,\rz)$ that
is commensurable with $\Sp(n,\gz)$,
the quotient $\hz_n/\Gamma$ carries a natural structure
as quasi projective algebraic variety. The Satake compactification
$\overline{\hz_n/\Gamma}$ is a projective variety which is closely related to
Siegel modular forms.
The Satake compactification  can be identified with the projective variety associated
to a graded algebra of modular forms.
We recall briefly its definition. A modular form $f$ of weight $r/2$, $r\in\gz$,
is a holomorphic function $f$ on $\hz_n$ with the transformation property
$$f(MZ)=v(M)\sqrt{\det(CZ+D)}^rf(Z)$$
for all $M\in\Gamma$.
In the case $n=1$ a regularity condition at the cusps has to be added.
Here $v(M)$ is system of complex numbers of absolute valued one, called the
multiplier system. It has to fulfil an obvious cocycle condition.
We denote this space by
$[\Gamma,r/2,v]$.
Fixing some starting weight $r_0$ and a multiplier system $v$ for it, we define the
ring
$$A(\Gamma):=\bigoplus_{r\in\gz}[\Gamma,rr_0/2,v^r].$$
This turns out to be a finitely generated graded algebra and its associated
projective variety $\proj A(\Gamma)$  can be identified with the Satake compactifictation.
The ring depends on the starting weight and the multiplier system but the assiciated
projective variety does not.
\smallskip
Basic examples of modular forms are given be theta series with respect to
characteristics.
By definition, a theta characteristic is an element $m={a\choose b}$ from
$(\gz/2\gz)^{2n}$. Here $a,b\in(\gz/2\gz)^n$
are column vectors. The characteristic is called even
if ${\transpose a}b=0$ and odd otherwise.
The group $\Sp(n,\gz/2\gz)$ acts on the set of characteristics by
$$M\{m\}:={\transpose M}\strut^{-1}m+\pmatrix{(A{\transpose B})_0\cr (C{\transpose D})_0}.$$
Here $S_0$ denotes the column built of the diagonal of a square matrix $S$.
It is well-known that $\Sp(n,\gz/2\gz)$ acts transitively on the subsets of
even and odd characteristics.
Recall that to any characteristic the theta function
$$\vartheta[m]=\sum_{g\in\gz^n}e^{\pii (Z[g+a/2]+{\transpose b}(g+a/2))}
\qquad (Z[g]=\transpose gZg)$$
can be defined.
Here we use the identification of $\gz/2\gz$ with the
subset $\{0,1\}\subset\gz$. It vanishes if and only if $m$ is odd.
Recall also that the formula
$$\vartheta[M\{m\}](MZ)=v(M,m)\sqrt{\det (CZ+D)}\vartheta[m](Z)$$
holds for $M\in\Gamma_n$, where $v(M,m)$ is a rather delicate
$8^{\hbox{\sevenrm th}}$ root of unity which depends on the choice of the
square root.
Sometimes we will use the notation
$$\vartheta[m]=\vartheta\Bigl[{a_1a_2\atop b_1b_2}\Bigr]
\quad\hbox{for}\quad m=\pmatrix{a_1\cr a_2\cr b_1\cr b_2}.$$
Following van Geemen and Nygaard, we consider the 8 functions
\medni
\halign{\qquad$\displaystyle#$\quad\hfil&$\displaystyle#$\quad\hfil
&$\displaystyle#$\quad\hfil&$\displaystyle#$\quad\hfil\cr
\vartheta\Bigl[{00\atop 00}\Bigr](Z),& \vartheta\Bigl[{00\atop 10}\Bigr](Z),&
\vartheta\Bigl[{00\atop 01}\Bigr](Z),& \vartheta\Bigl[{00\atop 11}\Bigr](Z),\cr
\noalign{\vskip3mm}
\vartheta\Bigl[{00\atop 00}\Bigr](2Z),& \vartheta\Bigl[{10\atop 00}\Bigr](2Z),&
\vartheta\Bigl[{01\atop 00}\Bigr](2Z),& \vartheta\Bigl[{11\atop 00}\Bigr](2Z).\cr}
\medni
If  we denote them by $Y_0,\dots, Y_3,X_0\dots,X_3$, then classical theta
relations show that the relations listed at the beginning of the introduction
hold.
The classical theta transformation formalism shows that the eight forms are modular
forms of weight $1/2$ for the group $\Gamma'$ introduced in the introduction and
that their multipliers on this group agree.
Since this is standard, we only give a short sketch of the proof.
The theta series $X_i$ are the so-called theta series of second kind.
One knows classically that they have the same multipliers
$\kappa(M)$ on the
group $\Gamma_2[2,4]$ (s.~for example [Ru]). By conjugation with the transformation
$Z\mapsto 2Z$ one shows that the $Y_i$ have the same multipliers
$\tilde\kappa(M)$ on the group $\Gamma_{2,0,\vartheta}[2]$.
We have to find the subgroup
$\Gamma'$ of $\Gamma_2[2,4]\cap \Gamma_{2,0,\vartheta}[2]$, where
$\kappa$ and $\tilde \kappa$ agree.
This just means that
$\vartheta[0](Z)\vartheta[0](2Z)$ and $\vartheta[0](2Z)^2$
have the same multipliers. But these are the standard thetas series with respect
to the binary forms $1\,0\choose 0\,2$ and $2\,0\choose 0\,2$.
The advantage of binary forms is that they have an even number of variables
and standard formula can be used, for example [Fr2], 7.1.\qed
\smallskip
We can express this in saying
that the $Y_0,\dots,X_3$ are contained in the ring
$$A(\Gamma'):=\bigoplus_{r\in\gz}[\Gamma',r/2,\kappa^r].$$
The precise result, slightly extending results of [GN], states:
\proclaim
{Proposition}
{Let be
$$\Gamma'=\{M\in\Gamma_2[2,4];\quad C\equiv 0\,\mod\,4,\
\hbox{\rm diag}\>C\equiv 0\,\mod\, 8,\ \det D\equiv\pm1\,\mod\,8\}.$$
The ring $A(\Gamma')$ is generated by the eight theta series above.
The relations
$$\eqalign{Y_0^2&=X_0^2+X_1^2+X_2^2+X_3^2,\cr
Y_1^2&=X_0^2-X_1^2+X_2^2-X_3^2,\cr
Y_2^2&=X_0^2+X_1^2-X_2^2-X_3^2,\cr
Y_3^2&=X_0^2-X_1^2-X_2^2+X_3^2\cr}$$
are defining relations. They define a subvariety $\calX$ of
$P^7(\cz)$ which can
be identified with the Satake compactification of $\hz_2/\Gamma'$.}
RingvGN%
\finishproclaim
{\it Proof.\/} Proposition 2.5 of [GN] says that this complete intersection is
the Satake compatification for some subgroup  of $\Sp(2,\gz)$.
Necessarily this must be a subgroup of what we called $\Gamma'$. An index
computation gives that they agree. The equality of the complete intersection
and the Satake compactification shows that $A(\Gamma')$ must be the
normalization of the factor ring $\cz[X_0,\dots,Y_3]$ by the
ideal generated by the above 4 relations. Using Serre's criterion for
normality, it follows that this ring is normal. This proves \RingvGN.\qed
\smallskip
In [GN] the modular form of weight 3
$$T=\vartheta\Bigl[{10\atop 00}\Bigr](Z)\vartheta\Bigl[{10\atop 01}\Bigr](Z)
\vartheta\Bigl[{01\atop 00}\Bigr](Z)\vartheta\Bigl[{01\atop 10}\Bigr](Z)
\vartheta\Bigl[{11\atop 00}\Bigr](Z)\vartheta\Bigl[{11\atop 11}\Bigr](Z)$$
has been introduced. The differential form
$$\omega=T\,dz_0\wedge dz_1\wedge dz_2,\qquad Z=\pmatrix{z_0&z_1\cr z_1&z_2},$$
is invariant under $\hat\Gamma_{2,0}[2]_{\Gotn}$
(The invariance under $\Gamma_{2,0}[2]_{\Gotn}$
has been proved in [FSM]. The behavior under the Fricke involution can be
derived from the following explicit formula.)
\proclaim
{Lemma}
{In terms of the coordinates $X_0,\dots X_3,Y_0\dots,Y_3$
we have that, up to a multiplicative constant, $\omega$ equals
$${X_2^4\over Y_0Y_1Y_2Y_3}d(X_0/X_2)
\wedge d(X_1/X_2)\wedge d(X_3/X_2).$$}
omegaXY%
\finishproclaim
{\it Proof.\/}
This is essentially the form of $\omega$, which has been derived in [FSM]
(see Theorem 4.5 and the formulae before it).\qed\smallskip
The zero locus of $T$ consists of the
fixed point sets of all $M\in \Gamma_{2,0}[2]_{\Gotn}$, which are conjugate inside
$\Sp(2,\gz)$ to the diagonal  matrix with diagonal
$(1,-1,1,-1)$. It can be checked that all these $M$ are contained in $\Gamma'$.
Hence we obtain the following well-known result:
\proclaim
{Remark}
{The differential form
$\omega$ defines a holomorphic differential form without zeros on
the smooth variety $\hz_2/\Gamma'$. It is invariant under
$\hat\Gamma_{2,0}[2]_{\Gotn}$.}
TCalDif%
\finishproclaim
\neupara{Automorphisms of the variety of van Geemen and Nygaard}%
As we will see (\calGgen), the group $\Gamma'$ is normal
in $\hat\Gamma_{2,0}[2]$. Hence this group acts
on the variety of van Geemen and Nygaard. We want to describe this action.
We will see that this action can be described by a linear action on the
variables $Y_0,\dots, X_3$, more precisely by a finite subgroup
of $\PGL(8,\cz)$. Recall that $\Gamma_{2,0}[2]$ is generated
by the matrices of the form
$$\pmatrix{E&S\cr0&E},\quad\pmatrix{U'&0\cr0&U^{-1}},\quad \pmatrix{E&0\cr2S&E}\qquad
(S=S'\ \hbox{integral}).$$
Let $M\in \hat\Gamma_{2,0}[2]$. For
$f\in[\Gamma',1/2,v_\vartheta]$ we set
$$f\vert M(Z)=\sqrt{\det(CZ+D)}^{-1/2} f(MZ).$$
The map $f\mapsto f\vert M$ is
an automorphism $\varphi_M$
of the 8-dimensional space spanned by $Y_0,\dots,X_3$.
It depends on the choice of a square root of $\det(CZ+D)$. Hence
$\pm\varphi_M$ is well defined.
Using standard theta transformation formulae we can compute these automorphisms
for the generators. It is sufficient to take the following 4.
\smallni
$$\eqalign{\hbox{matrix}\hbox to 1cm{\hfill}&\qquad
\hbox{corresponding transformation}\cr
\noalign{\vskip2mm}
%\pmatrix{E&S\cr0&E},\ S=\pmatrix{1&1\cr0&1}&\qquad
%(Y_1,Y_0,Y_3,Y_2,X_0,-\imag X_1,X_2,-\imag X_3)\cr
\pmatrix{\transpose U&0\cr0&U^{-1}},\ U=\pmatrix{1&1\cr0&1}&\qquad
(Y_0,Y_1,Y_3,Y_2,X_0,X_3,X_2,X_1)\cr
\pmatrix{^tU&0\cr0&U^{-1}},\ U=\pmatrix{1&0\cr1&1}&\qquad
(Y_0,Y_3,Y_2,Y_1,X_0,X_1,X_3,X_2)\cr
\pmatrix{E&0\cr S&E},\ S=\pmatrix{2&0\cr0&0}&\qquad
(Y_0,-\imag Y_1,Y_2,-\imag Y_3,X_1, X_0,X_3,X_2)\cr
\noalign{\vskip2mm}
J\quad \hbox{(Fricke involution)}
\hbox to10mm{\hfil}&\qquad
\hbox to 0pt{\hss$\sqrt2\cdot$}(X_0,X_1,X_2,X_3,
Y_0/2,Y_1/2,Y_2/2,Y_3/2)\cr
}$$
\proclaim
{Lemma}
{The group $\Gamma'$ is normal in $\Gamma_{2,0}[2]$.
The group $\calG$ generated by the transformations
$\pm\varphi_M$, $M\in\hat\Gamma_{2,0}[2]$, is already generated by the
above four transformations. It has order\/ $98\,304$.
The map $M\mapsto \pm\varphi_M$ defines a homomorphism
$$\hat\Gamma_{2,0}[2]\lo\calG/\pm.$$
The group $\calG$ contains the subgroup $\calZ$ of order
$4$ which is generated by multiplication with $\imag$.
The above homomorphism induces an isomorphism
$$\hat\Gamma_{2,0}[2]/\Gamma'\Isom\calG/\calZ\qquad
(\hbox{order }24\,576).$$
}
calGgen%
\finishproclaim
{\it Proof.\/} Since the generators of $\Gamma_{2,0}[2]$ act on $\calX$,
the group $\Gamma'$ must be a normal subgroup. The rest follows
by comparing indices.\qed
\neupara{The stabilizer of a node}%
The variety $\calX $ has 96 singularities, which all are ordinary
double points (nodes).
They are zero dimensional boundary points, but not each zero dimensional
boundary is singular. In coordinates the singularities
can be described as follows:
\smallskip
One node is given by
$$P=[\sqrt 2,0,\sqrt 2,0,\,1,1,0,0].$$
Changing signs it produces 8 nodes, which are characterized by the property
$Y_1=Y_3=X_2=X_3=0$. Similarly, one gets 8 nodes with
$Y_1=Y_3=X_0=X_1=0$. So one has 16 nodes with $Y_1=Y_3=0$.
In the same way one gets 16 nodes with the
property $Y_i=Y_j=0$ for each other pair $0\le i<j\le 3$.
This gives $96=8\cdot 16$ nodes. It is easy to check by hand that they exhaust
all singular points. This description of the nodes also shows:
\proclaim
{Lemma}
{The group  $\hat\Gamma_{2,0}[2]$ acts on
the\/ $96$ nodes is transitive.}
ActNode%
\finishproclaim
The following matrices
\medni
\halign{\quad\qquad$#$\qquad\hfil&$#$\hfil\cr
M_1=\pmatrix{1&0&0&0\cr0&1&0&0\cr0&2&1&0\cr2&0&0&1},&
%M_0=\pmatrix{1&0&0&0\cr0&1&0&0\cr4&0&1&0\cr0&0&0&1},\cr
%M_1=\pmatrix{1&0&0&0\cr0&1&0&2\cr0&0&1&0\cr0&0&0&1},&
M_2=\pmatrix{1&0&0&0\cr0&1&0&0\cr2&0&1&0\cr0&0&0&1},\cr
M_3=\pmatrix{1&0&0&0\cr0&1&0&1\cr0&0&1&0\cr0&0&0&1},&
M_4=\pmatrix{1&0&0&0\cr1&1&0&0\cr0&0&1&-1\cr0&0&0&1},\cr
M_5=\pmatrix{1&0&0&1\cr0&1&1&0\cr0&0&1&0\cr0&0&0&1},&
M_6=\pmatrix{0&1&0&0\cr1&0&0&0\cr0&0&0&1\cr0&0&1&0}\circ J
\cr
}
\medni
are contained in $\hat\Gamma_{2,0}[2]$ and stabilize the node $P$.
We consider the group, which is generated by them and $\Gamma'$.
One can check that the factor group mod $\Gamma'$ has order 128.
Together with \ActNode\ we obtain:
%\vskip0pt\RAND{A better argument of the above would be welcome}
\proclaim
{Proposition}
{The stabilizer $\hat\Gamma_{2,0}[2]_P$ of the node $P$
inside $\hat\Gamma_{2,0}[2]$ is generated by
$\Gamma'$ and the matrices $M_1,M_2,\dots,M_6$.
}
StabGen%
\finishproclaim
%\RAND{\StabGen\ needs comments on the proof}
In a neighborhood of $P$ we  can use the  affine coordinates
$$ \eta_0=Y_0/X_1,\, \eta_1=Y_1/X_1,\, \eta_2=Y_2/X_1,\, \eta_3=Y_3/
X_1,$$
$$ \xi_0=X_0/X_1,\,\ \xi_2=X_2/X_1,\, \xi_3=X_3X_1$$
Then substituting the  affine version of the third equation in the
fourth we get
$$\eta_1^2-\eta_3^2= 2(\xi_2^2-\xi_3^2)$$
Setting
$$x_1=\eta_1-\eta_3,\quad x_4=\eta_1+\eta_3,\quad
x_2=\sqrt 2(\xi_2 -\xi_3 ),\quad x_3=\sqrt 2(\xi_2 +\xi_3 ),$$
the relation gets the simple form
$$x_1x_4=x_2x_3.$$
So we have lead to
consider the quadric
$$Q:=\{(x_1,x_2,x_3,x_4)\in\cz^4;\quad x_1x_4=x_2x_3\}.$$
This is a three dimensional affine variety with a unique singularity at the
origin.
The above construction gives an \'etale map of germs
$$(\calX ,P)\lo (Q,0).$$
Sometimes we write the elements of $\cz^4$ as matrices
$$X=\pmatrix{x_1&x_2\cr x_3&x_4}.$$
Then $Q$ is defined by $\det X=0$.
The group $\GL(2,\cz)\times\GL(2,\cz)$ acts on the quadric by means of
$$X\loma AX\transpose B.$$
In this way we can consider $(\GL(2,\cz)\times\GL(2,\cz))/\cz^*$ as subgroup
of $\GL(4,\cz)$. Another transformation, which leaves the quadric invariant,
is $X\mapsto\transpose X$. We also can consider it as element of
$\GL(4,\cz)$.
\smallskip
We consider the transformations
\bigni
\halign{\quad$#$\hfil&\quad$#$\hfil\cr
m_1\pmatrix{x_1&x_2\cr x_3&x_4}=\pmatrix{x_4&x_2\cr x_3&x_1},&
%m_0\pmatrix{x_1&x_2\cr x_3&x_4}=\pmatrix{-x_1&x_2\cr x_3&-x_4}\cr
%m_1\pmatrix{x_1&x_2\cr x_3&x_4}=\pmatrix{x_1&-x_2\cr -x_3&x_4},&
m_2\pmatrix{x_1&x_2\cr x_3&x_4}=\pmatrix{-\imag x_1&-x_2\cr x_3&-\imag x_4},
\cr
m_3\pmatrix{x_1&x_2\cr x_3&x_4}=\pmatrix{-x_1&\imag x_2\cr \imag x_3&x_4},&
m_4\pmatrix{x_1&x_2\cr x_3&x_4}=\pmatrix{- x_1&-x_2\cr x_3&x_4},\cr
m_5\pmatrix{x_1&x_2\cr x_3&x_4}=\pmatrix{x_1&x_3\cr x_2&x_4},&
m_6\pmatrix{x_1&x_2\cr x_3&x_4}=
\pmatrix{x_2&x_1\cr x_4&x_3}.
\cr}%
\medni
They are contained in the extension of index two of the
image of the group
$(\GL(2,\cz)\times\GL(2,\cz))$,
which is generated by $X\mapsto{\transpose X}$.
Hence this subgroup acts on $Q$.
\proclaim
{Lemma}
{The transformations
$$m_1,m_2,m_3,m_4,m_5,m_6$$
generate a group $G$ of order\/ $256$.
}
DefGH%
\finishproclaim
Since the proofs of this Lemma and the
following Proposition can be given by computation, we omit  them.
\proclaim
{Proposition}
{The assignment
$$M_i\loma m_i\quad (1\le i\le 6)$$
induces an isomorphism
$$\hat\Gamma_{2,0}[2]_P/\Gamma'\Isom G.$$
The described identification of germs $(X(\Gamma'),P)$ and
$(Q,0)$ is equivariant.
}
IsoStabs%
\finishproclaim
We are interested in the subgroup of index two $\hat\Gamma_{2,0}[2]_{\Gotn}$.
We have to intersect it with the stabilizer.
It is easy to check that the elements
$$M_2^2,\ M_3^2,\ M_2M_1,\ M_3M_1,\ M_4M_1,\ M_5M_1,\ M_6M_1$$
are contained in $\hat\Gamma_{2,0}[2]_{\Gotn}$.
One also can check that
the elements
$$m_2^2,\ m_3^2,\ m_2m_1,\ m_3m_1,\ m_4m_1,\ m_5m_1,\ m_6m_1.$$
generate a group $H$ of order\/ $128$.
In this way one obtains:
\proclaim
{Proposition}
{The stabilizer of $P$ in $\hat\Gamma_{2,0}[2]_{\Gotn}$ is a subgroup of
index two of $\hat\Gamma_{2,0}[2]_P$. The restriction of \IsoStabs\ induces an isomorphism
$$(\hat\Gamma_{2,0}[2]_P\cap \hat\Gamma_{2,0}[2]_{\Gotn})/\Gamma'\Isom H.$$}
Isozstabs%
\finishproclaim
\neupara{Some results about crepant resolutions}%
In this section we recall some known result about projective crepant
resolutions.
\smallskip
We need the notion of a crepant resolution for certain
normal three dimensional varieties $X$.
\proclaim
{Definition}
{Let $X$ be a connected three dimensional normal complex space
with singular locus $S$.
Assume that for each point $a\in X$
there exists an open neighborhood $U$ and holomorphic three form $\alpha$
on $U-S$ without zeros.
A {\emph crepant resolution} $f:\tilde X\to X$ is a holomorphic map
of a 
connected (smooth) complex manifold $\tilde X$ onto $X$, such that
$\tilde X-f^{-1}(S)\to X-S$ is biholomorphic and such that
$\alpha$ extends to
a holomorphic three form without zeros on the inverse image
$\tilde U=f^{-1}(U)$.}
DefCrep%
\finishproclaim
The existence of a \it quasiprojective\/ \rm crepant desingularization is only a
local question  (in the three-dimensional
case):
\proclaim
{Lemma (Roan)}
{Under the assumptions of \DefCrep\ the following holds:
The existence of a  crepant
desingularization is granted, if
there exists an open covering $U_i\subset X$, such that each $U_i$ admits
a  crepant desingularization.}
LocCrep%
\finishproclaim
We reproduce the argument of Roan: The singular locus is a curve $S$.
Over the generic point of $S$ the crepant resolution is unique.
For this reason one can choose the crepant resolutions over the finitely many singular
points of $S$ arbitrarily and glue them to a global resolution.
\smallskip
We use
a result of a general result about the existence of a
resolutions of quotient singularities:
\proclaim
{Theorem}
{Let $X$ be a complex manifold of dimension
three and $G$ a finite group of automorphisms
of it. Assume that every point of $X/G$
admits an open neighborhood (in the analytic topology) such that
on its regular locus there exists a three-form without zeros.
Then $X/G$ admits a
crepant
desingularization.}
Ito%
\finishproclaim
Because of \LocCrep\ it is sufficient to prove this result in the sitation
$\cz^3/G$ where $G\subset\GL(3,\cz)$ is a finite subgroup. We can assume that
$G$ doesn't contain quasi reflections. Then the assumption in \Ito\
gives $G\subset\SL(3,\cz)$.
We refer to [Re] (especially section 5) for historical comments and
basic results.
\neupara{Some quotients of an ordinary double point}%
We will study  the node $(Q,0)$ and some of its quotients.
This is related to work of Davis [Da], where also several quotients
have been considered.
\proclaim
{Lemma}
{The restriction of
$$\alpha={1\over x_1^2-x_4^2}(x_1dx_1+x_4dx_4)\wedge dx_2\wedge dx_3$$
is a holomorphic differential form of degree three
on $Q-\{0\}$ without zeros.
If one identifies a small neighborhood of the origin in $Q$ with a small
neighborhood of $P\in X(\Gamma')$, we have $\omega=h\alpha$, where
$h$ is a holomorphic invertible function on this neighborhood.
}
Diffinv%
\finishproclaim
{\it Proof.\/}
We cover $Q$ by 4 charts corresponding to $x_i\ne 0$.
For example the part $x_1\ne 0$ is the  image of
$$\{(x_1,x_2,x_3)\in\cz^3;\ x_1\ne 0\}\lo\cz^4,\quad
(x_1,x_2,x_3)\loma (x_1,x_2,x_3,x_2x_3/x_1).$$
Pulling back $\alpha$ we get
$${1\over x_1^2-(x_2x_3/x_1)^2}(x_1dx_1
-x_2^2x_3^2/x_1^3 dx_1)\wedge dx_2\wedge dx_3={dx_1\wedge dx_2\wedge dx_3\over x_1}$$
This is holomorphic and without zeros on this chart. The other charts
are treated in a similar way.
\smallskip
Since $\alpha$ and $\omega$ both are 3-forms without zeros outside the singularity,
we get $\omega=h\alpha$, where $h$ is a holomorphic function without zeros
outside the singularity. Since isolated singularities of analytic
functions in more than one variable
cannot exist, $h$ and $h^{-1}$ are holomorphic also at the singularity.
\qed
\smallskip
The following result can be found in [Fri],
see also [Jo], 6.3.
\proclaim
{Proposition}
{The quadric $Q$ admits a small desingularization $\tilde Q\to Q$.
This means that $\tilde Q$ is  a smooth connected variety,
the inverse image of the node $0$ is a curve and the map from the
complement of this curve maps biholomorphically onto $Q-\{0\}$.
Such a desingularization is crepant.}
Mdes%
\finishproclaim
A small resolution is not unique.
Actually there exist two different isomorphy classes of such small resolutions.
They can be obtained by blowing up the ideals
$(x_1,x_3)$ or $(x_1,x_4)$ in $\cz[x_1,x_2,x_3,x_4]/(x_1x_4-x_2x_3)$.
From this explicit description one can derive:
\proclaim
{Lemma}
{The elements of the image of
$\GL(2,\cz)\times\GL(2,\cz)$ extend to biholomorphic
transformations of any small resolution, but the transformation
$$X\loma {\transpose X},$$
which is also an automorphism of $Q$, does not.}
Fort%
\finishproclaim

The group $G$ (see \IsoStabs)
is not contained in the image of $\GL(2,\cz)\times\GL(2,\cz)$.
The intersection with this group defines a subgroup
$G_0\subset G$ of index two.
It is generated by the elements $m_1m_5,m_2,m_3,m_4,m_6$.
One checks that these elements have determinant 1 (considered in $\GL(4,\cz)$).
Hence $G_0$ also can be defined as intersection of $G$ with $\SL(4,\cz)$.
We denote by $H_0$ the intersection of
$H$ and $G_0$. This is a group of order
$64$.
\proclaim
{Lemma}
{The groups $H$ and $H_0$ have the same center.
It is the group of order
$2$ generated by the
transformation
$x\mapsto-x$.}
Cent%
\finishproclaim
We omit the proof, since it can be done by simple computation.
\proclaim
{Lemma}
{The differential form $\alpha$ on $\cz^4$ (s.~\Diffinv)
is invariant under $H$.
Hence also the function $h$ in \Diffinv\ is $H$-invariant.}
Difinv%
\finishproclaim
{\it Proof.\/}
The invariance can be checked directly for the generators.\qed
\proclaim
{Theorem}
{Let $K$ be any subgroup of $H$. Then
the quotient $Q_0:=Q/K$ admits a {crepant desingularization}.}
Crep%
\finishproclaim
%\RAND{Proof must be checked carefully}
For the proof we have to differ between 4 types
of subgroups $K$.
\smallni 1) $K$ is contained in the subgroup $H_0$.\hfill\break
2) $K$ contains the transformation $x\mapsto -x$.
\vskip 1mm
\item{3)} $K$  contains one of the two transformations
$$\pmatrix{x_1&x_2\cr x_3&x_4}\loma \pmatrix{-x_1&x_2\cr x_3&-x_4}\quad
\hbox{or}\quad \pmatrix{x_1&-x_2\cr-x_3&x_4}$$
in its center.
\vskip1mm
\item{4)} $K$ is a the cyclic subgroup of order 4 that is
conjugated to one of the two following (given by a
generator of order 4):
$$\pmatrix{x_1&x_2\cr x_3&x_4}\loma
\pmatrix{x_2&x_4\cr x_1&x_3}\quad\hbox{or}\quad
\pmatrix{x_2&\imag x_4\cr -\imag x_1&x_3}.$$
\vskip0pt\noindent
These classes are not disjoint. But each subgroup
of $H$ is contained in at least one of the 4 classes.
This can be checked by hand or quicker by means of a computer.
\smallskip
We are going to discuss the 4 cases.
We begin with the
\smallni
{\it First type.\/} Because of \Fort\
in this case the action of $K$ extends to a small
resolution $\tilde Q\to Q$.
The differential form $\alpha$
extends to a holomorphic differential form without zeros on $\tilde Q$,
since singularities or zeros can only occur in codimension one
on a smooth variety.
By \Ito\ the variety $\tilde Q/K$ admits
a crepant resolution. Hence $Q/K$ also admits one.
\smallni{\it Second type.\/}
We start to blow up the origin of $\cz^4$. The group
$K$ still acts biholomorphically on this blow up. A typical chart
of the blow up is the $\cz^4$ with the coordinates
$$(t_1,t_2,t_3,x_4)=(x_1/x_4,x_2/x_4,x_3/x_4,x_4).$$
We consider in the blow up of $\cz^4$ the closed smooth subvariety
$\tilde Q$, which is defined in this chart by $t_1=t_2t_3$.
Its image in $\cz^4$ is $Q$. Hence $\tilde Q\to Q$ is just a
desingularization of $Q$. (Actually it is the blow-up of $Q$ at the
origin.) The chart of $\tilde Q$, which we consider, is a
$\cz^3$ with the coordinates $t_2,t_3,x_4$.
The differential form $\alpha$ in these coordinates can be computed.
Up to a constant factor it is
$$x_4\>dt_2\wedge dt_3\wedge dx_4.$$
So it gets a zero of order one along $x_4=0$.
The transformation $x\to -x$ just changes
the sign of each variable $x_i$. Hence it acts on
$t_2,t_3,x_4$  as reflection, which
changes the sign of the third variable only.
The quotient is a $\cz^3$ with the coordinates
$$(t_2,t_3,t_4)=(x_2/x_4,x_3/x_4,x_4^2).$$
Hence $\tilde Q/\pm$ is a smooth variety, the affine piece
in consideration
a $\cz^3$ with coordinates $t_2,t_3,t_4$. In these coordinates
$\alpha$ appears as holomorphic differential form without zeros.
By the general theorem \Ito, the quotient
$\tilde Q/K$ and hence $Q/K$ admits a crepant
resolution.
\smallni{\it Third type.\/}
The two cases are equivalent, hence we can assume that
$$\sigma(x_1,x_2,x_3,x_4)=(x_1,-x_2,-x_3,x_4)$$
is in the center of $K$.
The ideal $(x_2,x_3)$ is invariant under $K$, since it describes
the fixed point locus of $\sigma$ which is in the center of $K$.
Hence the action of $K$ extends
to an action by biholomorphic transformations on the blow up
$\calC$ of $\cz^4$ along this ideal. The manifold $\calC$ can be covered by
two $\cz^4$ using the coordinates $(x_1,x_2/x_3,x_3,x_4)$
and $(x_1,x_2,x_3/x_2,x_4)$.
We take the quotient of $\calC$ by $\sigma$.
The quotient $\calC/\sigma$ is covered by two $\cz^4$ with coordinates
$$(u_1,u_2,u_3,u_4)=(x_1,x_2/x_3,x_3^2,x_4)\quad\hbox{and}\quad
(v_1,v_2,v_3,v_4)=(x_1,x_2^2,x_3/x_2,x_4).$$
We consider the subvariety $Q'\subset \calC/\sigma$, which in the two
affine pieces is given by $u_1u_4=u_2u_3$ and $v_1v_4=v_2v_3$.
This variety has two singular points, which correspond to the origins
of the affine pieces and which are ordinary double points.
There is a natural projection $Q'\to Q/\sigma$. The group \cal$K$
acts on $Q'$.
We need the stabilizers of the two singularities.
\smallni
%\RAND{Basic point of the argument, must be checked}%
{\it Claim.\/} The stabilizer of the singularity $u_1=\dots=u_4=0$
acts by substitutions of the form
$$U\loma CU\transpose D,\quad\hbox{where}\quad U=\pmatrix{u_1&u_2\cr u_3&u_4}.$$
{\it Proof of the Claim.\/}
We consider an element of the stabilizer. It might be of the form
$$X\loma X\loma AX\transpose B\quad\hbox{or}
\quad X\loma A\transpose X\transpose B.$$
%\RAND{Basic point of the argument continued, must be checked}%
We have to use now that this transformation commutes with $\sigma$.
In each case
this means that $A$ and $B$ both are of the form
${*\,0\choose 0\>*}$ or both are of the form ${0\>*\choose *\,0}$.
Both cases are similar. For simplicity we take the case
$A={a\,0\choose 0\>d}$   and $B={\alpha\,0\choose 0\>\delta}$.
Then the transformation $AX\transpose B$ corresponds to
$CU\transpose D$, where
$$C=\pmatrix{1&0\cr0&d^2\alpha/a},\quad D=\pmatrix{a\alpha&0\cr
0&a\delta/d\alpha}.$$
But the transformation
$A\transpose X\transpose B$
interchanges the $u$- and $v$-chart and especially the two singularities
are interchanged. Hence this transformation is not contained
in the stabilizer.
\smallskip
We consider now the holomorphic map
$$Q'/\calK\lo Q/\calK,$$
which is induced by $Q'\to Q/\sigma$.
The differential form $\alpha$ can be considered as
meromorphic differential
form on $\calC/\sigma$. One checks that in the coordinates $u_1,u_2,u_3,u_4$
it is given by the same equation as the original $\alpha$, just replacing
the letters $x$ by $u$. The same is true for the coordinates $v_1,v_2,v_3,v_4$.
Hence $\alpha$ gives a holomorphic differential form without zeros
on the regular
locus of $Q'$.
The claim shows that $\calK$ extends to a suitable chosen small
resolution $\tilde Q'$ of $Q'$.
To be concrete one can blow up the irreducible surface,
which is defined
in the $u$-coordinates  by the ideal
$(u_2,u_3)$ and in the $v$-coordinates by $(v_2,v_3)$.
Then the group $K$ acts on $\tilde Q'$.
The differential form
$\alpha$ extends to a holomorphic differential form on
$\tilde Q'$ without zeros
and is
invariant under $\calK$.
\smallni
{\it Fourth type.\/}
We consider the first case,
$\sigma(x_1,x_2,x_3,x_4)=(x_2,x_4,x_1,x_3)$,
the second is similar.
It is better then to use the coordinates
$$y_1=x_1+x_4,\ y_2=x_1-x_4,\ y_3=x_2+x_3,\ y_4=x_2-x_3.$$
The quadric then takes the equation
$y_1^2+y_4^2=y_2^2+y_3^2$.
We blow up the ideal $(y_2,y_4)$.
We denote by $\calC$ the blow up of $\cz^4$ along this ideal.
One chart of the blow up is
$(y_1,y_2,y_3,y_4/y_2)$.
Taking quotient by $\sigma^2$ gives the chart
$$(u_1,u_2,u_3,u_4)=(y_1,y_2^2,y_3,y_4/y_2).$$
Now the the quadric
$y_1^2+y_4^2=y_2^2+y_3^2$
gets the form
$$u_1^2+u_2(u_4^2-1)-u_3^2=0.$$
This 3-fold has two singular points,
$(0,0,0,1)$ and $(0,0,0,-1)$.
Take the first one. For this point one can
take as local parameter
$$v_4:=u_4^2-1\qquad  (=(u_4-1)(u_4+1)).$$
Now the 3-fold is given by
$u_1^2+u_2v_4-u_3^2=0$. Hence the singularity is an
ordinary double point.
This consideration shows that the transform of the quadric appears
in $\calC/\sigma^2$ as 3-fold with two singular points, which are nodes.
The transformation $\sigma$ interchanges the two nodes and hence acts without
fixed points. The rest of the proof is analogously to the third case.
\smallskip
Now  the general result \Ito\ shows that
$\tilde Q'/\calK$ admits a crepant
resolution. This gives a crepant resolution of $Q/\calK$.\qed
\smallskip
We recall that we defined a local etale map
$$(X(\Gamma'),P)=(\calX,P)\lo (Q,0).$$
If $\Gamma$ is a group between $\Gamma'$ and $\hat\Gamma_{2,0}[2]_{\Gotn}$
and $\calK$ the corresponding subgroup of $\calH$, we still get a local etale
map
$$(X(\Gamma),P)\lo (Q/\calK,0).$$
Pulling back a crepant resolution of $Q/K$ we get a crepant resolution of
some affine neighborhood of $P$ in $X(\Gamma)$.
\smallskip
Now we can prove the a basic result of this paper, formulated in the
introduction, which
states that the Siegel threefold for each group between
$\Gamma'$ and
$\hat\Gamma_{2,0}[2]_{\Gotn}$ admits a (weak) Calabi-Yau model:
By \Ito\ there exists a crepant
resolution of the complement of the (images of the) nodes in
$X(\Gamma)$.
As we have just seen for each node there exists an open neighborhood which
admits a crepant resolution.
Hence we can apply \LocCrep\ to obtain a (not necessarily projective)
crepant resolution of $X(\Gamma)$.
\smallskip
It is a natural question whether a group $\Gamma$ extends to a group
of biholomorphic maps of a crepant resolution $\tilde \calX$. Such a resolution is
not unique. There exists a projective one but there also exist some which
are not projective. A necessary condition of $\Gamma$ to extend is that
the stabilizers of the nodes extend as described in \Fort. This is a condition which
can be checked. Assume that it is satisfied.
Then we can choose one node $a$ and desingularization of this
node. We can extend $\Gamma_a$ to this desingularization. Let now $g\in\Gamma$.
The choice of the resolution of $a$ dictates us the choice of the resolution at
$g(a)$ and the assumption about $G_a$ makes this choice independent of the
choice of $g$. In this way we resolve the whole orbit of $a$ and then in the
same way the other orbits. In this way we obtain:
\proclaim
{Lemma}
{Let $\Gamma$ be a group between $\Gamma'$ and
$\hat\Gamma_{2,0}[2]_{\Gotn}$. Assume that the stabilizer at an arbitrary node
satisfies the local condition \Fort. Then there exists a crepant resolution
$\tilde\calX\to\calX$ in the category of complex spaces such that $\Gamma$
extends biholomorphicallly to $\tilde X$.}
ExtCom%
\finishproclaim
As van Geemen pointed out to the authors it is usually not possible to
get $\tilde\calX$ as a projective variety. For example one can show that
there are freely acting groups of order $\#\Gamma/\Gamma'=32$.
For them it is not possible
to get a projective $\tilde\calX$.
\def\pic{{\rm pic}}
%x
\neupara{The Picard group}%
In the following we consider a group $\Gamma$ between
$\Gamma'$ and $\hat\Gamma_{2,0}[2]_{\Gotn}$.
We denote by $X(\Gamma)$ the Satake compactification of
$\hz_2/\Gamma$ and by $X(\Gamma)_{\hbox{\sevenrm reg}}$ its regular locus.
We want to compute the Picard number of this variety.
We denote it by
$$\pic(\Gamma):=\dim \Pic(X(\Gamma)_{\hbox{\sevenrm reg}})\otimes_\gz\qz.$$
Since the map $X(\Gamma)\to X(\Gamma')$ is unramified
in codimension one, we have
$$\Pic(X(\Gamma)_{\hbox{\sevenrm reg}})=
\Pic(X(\Gamma')_{\hbox{\sevenrm reg}})^\Gamma.$$
in [GN] the Picard number has been computed for $\Gamma'$:
$$\pic(\Gamma')=32.$$
The proof of van Geemen and Nygaard rests on
counting numbers of points with values in some finite fields and
making use of the Weil conjectures. Therefore it may be difficult
to get from these computations the action of $\hat\Gamma_{2,0}[2]$
on $\Pic(\calX _{\hbox{\sevenrm reg}})$.
For this reason we need some explicit description of
generating divisors.
\smallskip
We use Igusa's cusp form $\chi_{35}$ of weight 35
for the full Siegel modular group.
\proclaim
{Theorem}
{The group $\Pic(\calX _{\hbox{\sevenrm reg}})$
is generated by the irreducible components of the
zero divisor of $\chi_{35}$.}
picGN%
\finishproclaim
The proof needs some computer calculation, which rests on the
following informations about the structure of $\chi_{35}$.
First we recall that $\chi_{35}$ is the product
$$\chi_{35}=\chi_5\cdot\chi_{30}$$
of two forms of weight 5 and 30, also for the full modular
group but with respect to its non-trivial character.
The form $\chi_5$ can be defined as the product of the
ten theta constants. They can be produced as follows:
One starts with the most trivial theta constant
$$\vartheta[0](Z)=\sum_{g\in\gz^2}e^{\pii Z[g]}$$
and applies the full modular group to it. This gives 10
modular forms and there product is $\chi_5$ (up to a constant
factor).
The form $\chi_{30}$ can be constructed in a similar way.
\proclaim
{Lemma}
{If on applies the full modular group to $X_0=\vartheta[0](2Z)$,
on gets (up to constant factors)\/ $60$ modular forms
(living on $\Gamma_2[2,4]$ all with the same multipliers).
There product up to a constant is $\chi_{30}$.
Examples of forms in the orbit are
$$\vartheta[0](Z/2)\quad\hbox{and}\quad X_0,X_1,X_2,X_3.$$
The\/ $60$ modular forms can be written as linear combinations
of $X_0,\dots,X_3$.
}
chiDr%
\finishproclaim
The last statement is true, since the full modular group acts
on the space generated by $X_0,\dots,X_3$. This action has been
studied in detail by Runge [Ru].
\smallskip
So far we have seen that the zero locus of $\chi_{35}$
considered on $\calX $ splits into the sum of 70 pairwise different
divisors. But these 70 divisors need not to be irreducible.
\smallskip
Now computer algebra comes into the game. Since we know the equations
of the 70 divisors in $P^7(\cz)$ we can  decompose them into
irreducibles by using the facility of computer algebra
to compute the primary decomposition of an ideal.
We have to be a little careful, since computer algebra works well
not over $\cz$ but only over a finitely generated field.
Hence we have to use a number field $K$.
We use $K=\qz(\zeta_8)$, where $\zeta_8$ is a primitive 8$^{\hbox{\sevenrm th}}$
root of unity. We got the following result using {\dunh MAGMA} [BMP]:
\proclaim
{Proposition}
{We consider $\calX $ as variety over $K=\qz(\zeta_8)$ (using the equations
of van Geemen and Nygaard). The zero locus of $\chi_{35}$ is defined
over this field. It splits over $K$ into\/ $132$ irreducible
components. More precisely we have:
\smallni
{\rm 1)} The theta constants $Y_0,Y_1,Y_2,Y_3$ have irreducible
divisors. The other\/ $6$ decompose into two irreducibles.
Hence $\chi_5$ contributes with\/ $16=4+2\cdot 6$ irreducible components.
\smallni
{\rm 2)} The forms $X_0,X_1,X_2,X_3$ have irreducible divisors.
The other\/ $56$ factors of $\chi_{30}$ decompose into pairs
of irreducibles. Hence $\chi_{30}$ contributes with\/ $116=4+2\cdot 56$
irreducible components.}
decIrr%
\finishproclaim
We conjecture that these 132 components are irreducible over
$\cz$. But there is no need for us to check this.
\smallskip
The proof of \picGN\ now can be given as follows.
Using Poincar\`e duality it is sufficient to construct a system
of curves $C_1,\dots,C_m$, which are complete and contained in
the regular locus of $\calX $ (i.e.\ they don't contain one
of the 96 nodes) and such that the intersection matrix between the
132 divisors and these curves has rank 32. The construction
of these curves can be given (again by means of a computer) as
follows: Take pairwise intersections of the 132 divisors, decompose
them into irreducibles and single out those components, which don't
contain  nodes.
\smallskip
In this way \picGN\ can be proved.
This explicit description of the Picard group allows us to
describe the action of the group $\hat\Gamma_{2,0}[2]$
on it. The group $\calG$ acts in a natural way on the ring
$\cz[Y_0,\dots,Y_3,X_0,\dots,X_4]$ and on its ideals.
Hence we can describe the action of $\calG$ on
$\Pic(\calX _{\hbox{\sevenrm reg}})$ explicitly.
Of course $\calZ$ acts trivially. Using the ismomorphism
$\hat\Gamma_{2,0}[2]/\Gamma'\cong\calG/\calZ$ we get the action
of $\hat\Gamma_{2,0}[2]$ on $\Pic(\calX _{\hbox{\sevenrm reg}})$.
Using the character table for $\calG/\calZ$ which can be produced
by computer algebra we get the decomposition into irreducibles.
\proclaim
{Theorem}
{The space $\Pic(\calX _{\hbox{\sevenrm reg}})\otimes_\gz\qz$
decomposes under $\hat\Gamma_{2,0}[2]$ into four irreducible
components of dimensions $1,3,12,16$.}
DecIrr%
\finishproclaim
These numbers can be explained as follows:
\vskip1mm\item{1)}
The 1-dimensional component comes from the  divisor
of a modular form.
\vskip1mm\item{2)}
The 3-dimensional space comes from the components of the 6
theta constants, which are different form $Y_0,\dots,Y_3$.
\vskip1mm\item{3)}
The 56 factors of $\chi_{30}$, which are different from $X_0,\dots,X_3$
decompose under $\Gamma_{2,0}[2]$ into two orbits of 24 and 32 elements.
Their irreducible components produce the spaces of dimension 12 and
16.
\smallni
This explicit picture of the action of $\hat\Gamma_{2,0}[2]$ on
$\Pic(\calX _{\hbox{\sevenrm reg}})$ allows to compute the
number $\pic(\Gamma)$ for every group $\Gamma$ between
$\Gamma'$ and $\hat\Gamma_{2,0}[2]$ and this can be done by means of a
program.
\neupara{Involutions}%
As we have seen, the group $\hat\Gamma_{2,0}[2]/\Gamma'\cong\calG/\calZ$
has order $24\,576$. We are interested in its subgroups of order
two. One can compute that there are 18 conjugacy classes of such subgroups,
and one can show that 10 of them are in the image of
$\hat\Gamma_{2,0}[2]_{\Gotn}$. In the following we list them.
There are two possibilities to define such a group. We could describe it by
an element $M\in\hat\Gamma_{2,0}[2]$,
such that the image of $M$ in $\hat\Gamma_{2,0}[2]/\Gamma'$
generates the group.
In the case $M\in\Gamma_{2,0}[2]$ it is enough to consider its image
in $\Sp(2,\gz/8\gz)$,
since $\Gamma'$ contains $\Gamma_2[8]$. The other possibility is to
give a matrix $g\in\calG$ such that its image in $\calG/\calZ$ generates
the subgroup of $\calG/\calZ$. Such a $g$ is determined up to a
power of $\imag$.
We give the 10 groups the names $\calG2\_\imag$, $1\le i\le 10$.
\bigni
\medni\halign{Group $\calG2\_#$\quad\hfil&$#$\quad\qquad\hfil&$#$\cr
1&
\pmatrix{3&0&4&0\cr
0&1&0&0\cr
0&0&3&0\cr
0&0&0&1}%
&
(Y_0,Y_1,Y_2,Y_3,-X_0,-X_1,-X_2,-X_3)\cr
\noalign{\vskip2mm}
2&
\pmatrix{5&2&6&2\cr
2&1&2&6\cr
4&4&1&6\cr
4&4&6&5}
&
%(Y_0,Y_1,Y_2,Y_3,X_0,X_1,X_2,X_3)\cr
(Y_0,-Y_1,-Y_2,Y_3,X_0,-X_1,-X_2,X_3)\cr
\noalign{\vskip2mm}
3&
\pmatrix{1&0&2&6\cr
2&1&2&6\cr
0&0&1&6\cr
0&0&0&1}
&
(Y_0,Y_1,Y_2,Y_3,X_0,-X_1,-X_2,X_3)\cr
\noalign{\vskip2mm}
4&
\pmatrix{3&6&4&2\cr
4&7&6&2\cr
0&0&3&4\cr
0&4&2&7}
&
(-Y_0,-Y_1,Y_2,Y_3,X_0,X_1,-X_2,-X_3)\cr
\noalign{\vskip2mm}
5&
\pmatrix{3&2&6&7\cr
2&3&7&2\cr
4&2&1&2\cr
2&4&2&1}%
&
(-Y_0,-Y_1,-Y_2,Y_3,X_0,X_1,X_2,-X_3)\cr
\noalign{\vskip2mm}
6&
\pmatrix{5&2&6&7\cr
0&3&1&0\cr
0&2&3&0\cr
6&4&6&5}
&
(-Y_0,-Y_1,-Y_2,Y_3,-X_0,-X_1,-X_2,X_3)\cr
\noalign{\vskip2mm}
7&
\pmatrix{7&4&7&6\cr
0&7&2&3\cr
6&4&5&0\cr
4&6&4&5}
&
(-Y_3,-\imag Y_2,-\imag Y_1,Y_0,X_3,\imag X_2,\imag X_1,-X_0)\cr
\noalign{\vskip2mm}
8&
\pmatrix{3&7&5&2\cr
6&1&2&6\cr
0&6&1&6\cr
2&0&7&3}
&
(-Y_1,Y_0,-\imag Y_2,-\imag Y_3,X_2,-\imag X_1,-X_0,\imag X_3)\cr
\noalign{\vskip2mm}
9&
\pmatrix{1&0&7&0\cr
0&7&0&0\cr
2&0&7&0\cr
0&0&0&7}
&
(-\imag Y_1,Y_0,-\imag Y_3,Y_2,X_1,-\imag X_0,X_3,-\imag X_2)\cr
\noalign{\vskip2mm}
10& \hbox{Fricke involution}&
\hbox to 0pt{\hss$\sqrt2\cdot$}(X_0,X_1,X_2,X_3,
Y_0/2,Y_1/2,Y_2/2,Y_3/2)
\cr
}
%In Termen der Programme fsrdach entsprechen werden die Gruppen
%G2_i erzeugt von Elementen sigma_j in der Riehnfolge
%1,2,4,3,6,5,7,?,?,10,  also z.B. G2_3 erzeugt von sigma4
\smallni
\proclaim
{Proposition}
{The following table
gives the numbers\/ $\dim\Pic(X(\Gamma)_{\hbox{\sevenrm reg}})\otimes_\gz\qz$
for the groups $\Gamma$ corresponding to the groups
$\calG2\_i$, $1\le i\le 10$
and the dimension of the fixed point locus of the generating involution
in each case.
(Dimension -1 means that the locus is empty.)
}
somPp%
\finishproclaim
\vbox{
\halign{#\qquad\qquad\hfil&\rm#\ \hfil
&\rm#\ \hfil&\rm#\ \hfil&\rm#\ \hfil&\rm#\
&\rm#\ \hfil&\rm#\ \hfil&\rm#\ \hfil&\rm#\ \hfil&\rm#\ \hfil&\rm#\cr
Picard number&16&24&16&16&16&20&20&16&16&18\cr
Dimension&-1&0&1&-1&-1&1&1&1&1&1\cr}}
\medni
We want to compute the Picard- and Euler number for a
crepant resolution.
The numbers of a small resolution $\tilde\calX$ of $\calX$ are
[GN], [CM]
$$e=64,\quad \pic=32.$$
Hence $\tilde\calX$ is a rigid Calabi-Yau manifold.
\proclaim
{Lemma}
{The action of the groups $\calG2\_i$ extends to a suitable
crepant desingularization $\tilde\calX$.}
ActExt%
\finishproclaim
The proof rests on Lemma \ExtCom. We omit details.
\qed
\smallskip
Each of the groups $\calG2\_i$ is generated by an involution $\sigma_i$.
We need information by the fixed point locus.
It can be checked that there is no component of dimension two.
We also know that the Calabi-Yau 3-form is invariant under $\sigma_i$.
This implies that the action of $\sigma_i$ on the tangent space
of a fixed point can be diagonalized with diagonal $(-1,-1,1)$.
From this follows:
\proclaim
{Lemma}
{The fixed point locus of $\sigma_i$ on $\tilde\calX$ is the
disjoint union of
smooth curves.
They are in one to one correspondence with the irreducible
components of the fixed point locus on $\calX$.
The fixed point locus on $\calX$ consists of curves and
isolated points which are nodes. If a node occurs
as isolated fixed point, then
the exceptional line over it is in the fixed point locus
on $\tilde\calX$.
}
FixCur%
\finishproclaim
The image of the fixed point locus of $\sigma_i$ in
$\tilde\calX/\sigma_i$ is curve, which is in the singular locus.
We claim that in a crepant resolution over each of its
components
there is only one exceptional divisor.
Locally around a fixed point the involution can be
described by
$(w_1,w_1,w_3)\mapsto (-w_1,-w_2,w_3)$. The crepant resolution
of the quotient of $\cz^3$ by this involution is easy
to describe (we did it in [FSM]) and one sees from this description
that the exceptional divisor is smooth and connected.
\proclaim
{Lemma}
{The Picard number of a crepant resolution of $\calX/\calG2\_i$ equals
the sum of the Picard number of the regular locus of $\calX/\calG2\_i$
(see \somPp) and the number of irreducible components of
the fixed point locus of $\sigma_i$ (considered in $\calX$ is enough).}
PicNumb%
\finishproclaim
We also have to compute the Euler number of a crepant resolution
of $\calX/\calG2\_i$. This is given by the string theoretic
Euler number $e(\tilde\calX ,\calG2\_i)$. We refer to [Re] for some
comments about this.
We recall the definition of $e(M,G)$. Here $G$ is a finite group
acting on a compact differentiable manifold $M$.
One has to consider the subset of $G\times G$ of all commuting pairs
$(g,h)$. Then the string theoretic Euler number is defined as
$$e(M,G)={1\over\#G}\sum_{gh=hg}e(M^{\langle g,h\rangle}).$$
Here $M^{\langle g,h\rangle}$ denotes the common fixed point set
of $g,h$.
The string theoretic Euler number has the following basic
property:
Assume that $M$ is a Calabi  Yau manifold and that $G$ acts by
biholomorphic transformations, which leave the
Calabi-Yau 3-form invariant. Assume that for $a\in M$ the stabilizer
$G_a$ acts on the tangent space as subgroup of the special linear group.
Then there exists a crepant desingularization of $M/G$
and for each
such desingularization its usual Euler number is $e(M,G)$.
\smallskip
We apply this formula in the case, where the order of $G$ is two.
There are 4 commuting pairs
$(e,e),\ (\sigma,e),\ (e,\sigma),\ (\sigma,\sigma)$.
\proclaim
{Lemma}
{The
Euler number of a crepant resolution of
$\tilde\calX/\calG2\_i$ is
$$e=32+{3\over 2}\sum_Ce(C),$$
where $C$ runs through the components of the
fixed point locus of $\sigma_i$ in $\tilde\calX$.}
EulBer%
\finishproclaim
The fixed point sets are easy to determine.
The involution can be considered as a linear transformation
$A:\cz^8\to\cz^8$, where $A^2=aE$. We want to consider the fixed
point locus of $A$ on $P^7(\cz)$. It corresponds to the
eigenspaces of $A$. The possible eigenvalues are the two
square roots of $a$. We denote the two eigenspaces by $V^+$ and $V^-$.
Hence $\cz^8=V^+\oplus V^-$. The projective spaces $P(V^+)$ and
$P(V^-)$ are two disjoint linear subspaces of $P^7(\cz)=P(\cz^8)$.
To get the fixed point set of $A$ inside the variety $\calX$ we have to
intersect this variety with the two linear subspaces.
Hence the fixed point set inside $\calX$ is the disjoint
union of two parts, where each of the parts can be empty of course.
\smallskip
Following these lines one gets:
\proclaim
{Proposition}
{The following table
describes fixed point sets of the involutions $\sigma_i$, $i\le i\le 11$,
on $\calX$ and the Hodge numbers of a Calabi-Yau model
of the quotient $\calX/\calG2\_i$.
\smallni\vbox{
\halign{\qquad\qquad$\sigma_{#}$:\quad\hfil&\hbox{\rm #}\qquad\hfil
&$#$\quad\hfil&$#$\hfil\cr
\omit&fixed points&\pic&e\cr
\noalign{\vskip1mm}
1& empty set&16&32\cr
2& 16 nodes&40&80\cr
3& 4 elliptic curves&20&32\cr
4& empty set&16&32\cr
5& empty set&16&32\cr
6& 8 conics in  planes ($\cong P^1$)&28&56\cr
7& 8 lines ($\cong P^1$) &28&56\cr
8& 2 elliptic curves&18&32\cr
9& 2 elliptic curves&18&32\cr
10&4 conics in planes ($\cong P^1$)&22&44\cr
}}
}
FixSet%
\finishproclaim
The equations for the fixed point loci can be given explicitly.
We just give as an example the 4 elliptic curves which describe
the fixed point locus of $\sigma_3$: They are described by the following
4 ideals:
$$\eqalign{
        &(Y_0 + Y_1,\quad Y_2 + Y_3,\quad X_1,\quad X_3,\quad
        Y_1^2 - X_0^2 - X_2^2,\quad
        Y_3^2 - X_0^2 + X_2^2),\cr
        &(Y_0 + Y_1,\quad Y_2 - Y_3,\quad X_1,\quad X_3,\quad
        Y_1^2 - X_0^2 - X_2^2,\quad Y_3^2 - X_0^2 + X_2^2),\cr
        &(Y_0 - Y_1,\quad Y_2+Y_3,\quad X_1,\quad X_3,\quad
        Y_1^2 - X_0^2 - X_2^2,\quad
        Y_3^2 - X_0^2 + X_2^2),\cr
        &(Y_0 - Y_1,\quad Y_2 - Y_3,\quad X_1,\quad  X_3,\quad
        Y_1^2 - X_0^2 - X_2^2,\quad
        Y_3^2 - X_0^2 + X_2^2).\cr}$$
\vfill\eject\noindent
{\paragratit References}%
\bigni
\item{[BH]} Borisov, L.\ Hua, Z.:
{\it On Calabi-Yau threefolds with large nonabelian fundamental groups,\/}
Proc. Amer. Math. Soc. {\bf 136}, 1549--1551 (2008)
\medskip
\item{[BMP]} Bosma,\ W., Cannon,\ J., Playoust,\ C.:
{\it The Magma algebra system. I. The user language,\/}
J. Symbolic Comput. {\bf24} (3-4), 235--265 (1997)
\medni
\item{[CM]} Cynk, S., Meyer, C.:
{\it Modular Calabi-Yau Threefolds of level eight,\/}
Internat.~J.~Math.\ {\bf 18}, no.~3, 331--347 (2007)
\medskip
\item{[CFS]} Cynk, S., Freitag, E., Salvati-Manni, R.:
{\it The geometry and arithmetic of a
Calabi-Yau  Siegel threefold,\/} preprint 2010
\medskip
\item{[Da]} Davis, R.:
{\it Quotients of the conifold in compact Calabi-Yau threefolds,
and new topological transitions,\/}
eprint arXiv: 0911.0708 [hep-th] (2009)
\medskip
\item{[Fr1]} Freitag, E.: {\it Siegelsche Modulfunktionen,} Grundlehren
der mathematischen Wissenschaften, Bd. {\bf 254}. Berlin Heidelberg New
York: Springer (1983)
\medskip
\item{[Fr2]} Freitag, E.:
{\it Singular modular forms and theta relations,\/}
Lecture notes in Math.\ {\bf 1487},
Springer-Verlag,
Berlin Heidelberg New
York (1991)
\medskip
\item{[Fri]} Friedmann, R.:
{\it Simultaneous resolution of threefold double points,\/}
Mathematische Ann.\ {\bf 274}, 671--689 (1986)
\medskip
\item{[FSM]} Freitag, E. Salvati Manni, R.:
{\it Some Siegel threefolds with a Calabi-Yau model,\/} preprint 2009
\medskip
\item{[GN]} van Geemen, B., Nygaard, N.O.:
{\it On the geometry and arithmetic of some Siegel modular threefolds,\/}
Journal of Number Theory {\bf 53}, 45--87 (1995)
\medskip
\item{[Ig]} Igusa, I.: {\it A desingularization problem in the theory of
Siegel modular functions,\/}  Math. Annalen {\bf 168} 228--260
(1967)
%\item{[Li]} Lin, H.W.:
%{\it On crepant resolution of some hypersurface singularities
%and a criterion for UFD,\/}
%Transactions of the Am.\ Math.\ Soc. Vol.\ {\bf 354}, No.~5,
%1861--1868 (2002)
%\medskip
\item{[Re]} Reid, M.: {\it La correspondence de McKay,\/}
S\'eminaire Bourbaki 1999/2000. Ast\'erisque No. {\bf 276}, 53--72 (2002)
\medskip
%\medskip
%\item{[Ru2]} Runge, B.:
%{\it On Siegel modular forms, part II\/}, Nagoya Math. J. {\bf138}, 179-197 (1995)
\item{[SM]} Salvati Manni, R.:
{\it Thetanullwerte and stable modular forms,\/}
Am. J. of Math. {\bf 111}, 435--455 (1989)
\medskip
%\item{[We2]} Weissauer, R.: {\it The Picard group of Siegel modular threefolds,\/}
\bye

\neupara{Some more cases}%
We consider the unimodular substitution $\tau$.
$$Z\mapsto\transpose UZU,\qquad U=\pmatrix{2&1\cr 1&1}.$$
It is contained in $\hat\Gamma_{2,0}[2]_{\Gotn}$.
Its image in $\hat\Gamma_{2,0}[2]_{\Gotn}/\Gamma'$ has  order three.
The action on the variables $Y_0,\dots,X_3$ is given by
$$\tau_3:\qquad(Y_0,Y_3,Y_1,Y_2,X_0,X_2,X_3,X_1).$$
It can be shown that each subgroup of order three is conjugated
inside $\hat\Gamma_{2,0}[2]/\Gamma'$ to this group.
We compute the Hodge numbers for a crepant resolution
of $\calX/\tau_3$.
The Picard number of the regular locus
can be computed by means of the results of section 6.
The result is 12.
It can be checked that $\tau_3$ extends to any small resolution
$\tilde X$. The fixed point set of $\tau_3$ in $\calX$ has 5 components.
One component is an elliptic curve the other 4 are isolated points
in the regular locus of $\calX$.
The elliptic curve looks locally like $\cz^3/\tau$, where $\tau$
acts by
$$\tau(w_1,w_2,w_3)=(\zeta w_1,\zeta^{-1}w_2,w_3),\quad \zeta^3=1.$$
The 4 isolated singularities
of $\calX/\tau$ look locally like $\cz^3/\tau$, where now $\tau_3$
acts by multiplication with a third root of unity.
It is easy to check that in a crepant resolution of this singularity
one exceptional divisor occurs. So we obtain
the Picard number $12+5=17$ for a crepant resolution of $\tilde\calX/\tau$.
The Euler number can be computed by means of the string theoretic formula:
There are nine commuting pairs. The pair $(e,e)$ gives the contribution
$64/3$ and the 8 other pairs give the contribution
$(8/3)e(\Fix(\tau_3))$, where $\Fix(\tau_3)$ is the fixed point
locus of $\tau_3$ in $\tilde\calX$.
But since this fixed point locus is smooth it also consists of one elliptic
curve and 4 isolated points. Hence we get $e(\Fix(\tau_3))=4$.
Collecting the data we get:
\proclaim
{Theorem}
{Let $\Gamma$ be the group generated by $\Gamma'$ and
the unimodular transformation
$$Z\mapsto\transpose UZU,\qquad U=\pmatrix{2&1\cr 1&1}.$$
This a normal extension of index three. The Hodge numbers of a
crepant resolution of $X(\Gamma)$ are
$$h^{11}=18,\quad h^{12}=2.$$
So we obtain a non rigid Calabi Yau variety.}
orddrei%
\finishproclaim
This example admits a refinement We consider the transformation
$$\tau_6:\qquad(\imag Y_0,\imag Y_4,\imag Y_1,\imag Y_2,
-\imag X_0,-\imag X_2,-\imag X_3,-\imag X_1)$$
on $P^7(\cz)$. It has order 6.
It can be shown that it is the image of a transformation
in $\hat\Gamma_{2,0}[2]_{\Gotn}$ and that it extends to
$\tilde\calX$. Actually its square generates the group of order
3 and its third power generates
$\calG2\_1$. The Picard number of the regular locus of $\calX/\tau_6$ computes
as $6$. The point now is that only $\tau_6^k$ for even $k$ has a non
empty fixed point set. Hence we can use the computation for $\tau_3$.
In this way we obtain $h^{11}=6+5=11$ for the Picard number of a crepant
resolution of $\calX/\tau_6$. The string theoretic formula for the
Euler number gives
$64/6+(8/6)4=16$. Hence $h^{12}=3$.
\proclaim
{Theorem}
{The transformation
$$\tau_6:\qquad(\imag Y_0,\imag Y_4,\imag Y_1,\imag Y_2,
-\imag X_0,-\imag X_2,-\imag X_3,-\imag X_1)$$
is in the image of
$\hat\Gamma_{2,0}[2]_{\Gotn}$. It generates a group of order 6 acting
on $\calX$. The quotient $\calX/\tau_6$ admits a crepant resolution
with Hodge numbers
$$h^{11}=11,\quad h^{12}=3.$$
Hence we obtain a non rigid Calabi Yau manifold.}
ord6%
\finishproclaim

\neupara{A particular case}%
We have to consider translation matrices
$$T_S=\pmatrix{E&S\cr0&E}$$
of level two, $S\equiv 0$ mod 2. Such a translation matrix is called
{\it reflective\/} if $S$ is congruent 0 mod 4 to one of the tree
$$\pmatrix{2&0\cr0&0},\quad \pmatrix{0&0\cr0&2},\quad \pmatrix{2&2\cr2&2}.$$
Actually reflective translations act as reflections on the Igusa desingularization
of level four,
$\tilde X(4)$, which explains the notation.
\proclaim
{Lemma}
{
The group
$\Gamma=\Gamma_2[2]\cap\Gamma_{2,0}[4]$ is generated by
\vskip1mm
\item{\rm 1)} The group $\Gamma_2[4]$,
\item{\rm 2)} The elements of $\Gamma_{2,0}[2]_{\Gotn}$, which are
conjugate inside $\Gamma_2$ to the diagonal matrix with diagonal
$(1,-1,1,-1)$.
\item{\rm 3)} All elements of $\Gamma_{2,0}[2]_{\Gotn}$, which are
conjugate inside $\Gamma_2$ to a reflective
translation matrix ${E\,S\choose0\,E}$ of
$\Gamma_2[2]$.
\vskip0pt
}
DGa%
\finishproclaim
The proof  can be easily done with the help of a computer.
The group $\Sp(4,\gz/4\gz)$ can easily implemented in the computer
algebra system {\dunh MAGMA} as matrix groups generated by the standard generators.
The subgroup $\Gamma_{2,0}[2]$ can be defined as subgroup by
standard generators. Using the command ``LowIndexSubgroups'' the subgroup
of index two of $\Sp(4,\gz/4\gz)$ can be implemented and the non trivial character
can be defined as homomorphism of this group into a cyclic group of order
two. Now it is possible do define the character $\chi_{\Gotn}$ as homomorphism
of $\Gamma_{2,0}[2]$ into the cyclic group of order two. The command ``Kernel'' gives
the kernel of the group. Then one simply runs through all elements of this
kernel and asks whether it is conjugated to one of the four elements occurring in
2) and 3). Then one takes the subgroup, which is generated by these elements and
verifies that it equals the image of $\Gamma_2[2]\cap\Gamma_{2,0}[4]$.\qed
\smallskip
The lemma  is similar to lemma 1.4 in [FS]. There the group
$\Gamma_{2,0}[2]_{\Gotn}\cap\Gamma_2[2]$ has been characterized by the
same properties 1)--3) with the only difference that the word
``reflective'' has been skipped. The same proof as in [FS] works with this weaker
assumption and gives the result, which we formulated in the introduction, namely
that the quotient of the Igusa desingularization for the principal congruence subgroup
of level four
$\tilde X(4)/(\Gamma_2[2]\cap\Gamma_{2,0}[4])$ admits a desingularization, which is
a Calabi-Yau manifold. The same then is true for any group between
$\Gamma_2[2]\cap\Gamma_{2,0}[4]$ and $\Gamma_{2,0}[2]_{\Gotn}$.
We give the structure of this distinguished example:
\proclaim
{Proposition}
{The  ring of modular forms
of even weight (with trivial multipliers)
for the group $\Gamma_0[4]\cap\Gamma[2]$
is generated by
$$ \tt{00}{00} \tt{00}{01}  \tt{00}{10} \tt{00}{11},$$
all pairs  of the form
$$ \tt{00}{ab}^2 \tt{00}{cd}^2$$
and
the $\vartheta_m^4$.
\smallni
If one want the generators also in the  odd weights,
it is enough to add the form  $T$ of weight 3.}
ReMf%
\finishproclaim
To simplify the equations we consider the ring of forms of even weights:
\proclaim
{Proposition}
{The  ring $A(\Gamma_0[4]\cap\Gamma[2] )^{(2)}$ in the even weights
 is equal to
$$ \cz\Bigl[ \tt{00}{00} \tt{00}{01}  \tt{00}{10} \tt{00}{11},  \tt{00}{
00} ^2,  \tt{00}{01}^2,   \tt{00}{10}^2,  \tt{00}{11}^2, y_4 \Bigr]^{(2)}$$
with
$$ y_4=-\tt{10}{
01}^4 -\tt{00}{
11}^4$$}
Rmz%
\finishproclaim
Denoting the above variables with  $y_5, x_0, x_1, x_2, x_3$ we have the ring
$$ \cz[y_5,  x_0, x_1, x_2, x_3 , y_4 ]^{(2)}$$
with  $x_i$ of weight 1 and $y_j$ of weight 2
We have  also the
 following defining relations
$$
 y_5^2= x_0x_1x_2x_3,
 $$
$$
 2y_5^2= x_0^2x_1^2+x_0^2x_3^2+ x_1^2x_3^2+  ( -x_2^2+x_0^2+x_1^2+x_3^2+y_4)y_4.
$$
We want to compute the Hodge numbers of the Calabi-Yau model.
\smallskip
We need some information about the group
$K:=\Gamma/\Gamma'$.
The basic information is that $K$ is abelian of order 32 and that
all elements are of order two.
So their fixed point loci are known form the previous section.
We know that they all
extend to a small resolution $\tilde\calX$.
We also know from the previous
section that the fixed point locus is curve $C\subset\tilde\calX$.
The image of $C$ in $\tilde\calX/K$ is the singular locus.
The local structure of a singularity is of the type
$\cz^3/A$, where $A$ either is a group of order 2, generated by
a transformation, which changes two signs or the group of order 4 which
cointains all sign changes at two positions. It is easy to describe
the crepant resolution for these singularities (see [FSM]) and from
this description on can see:
\proclaim
{Lemma}
{The number of exceptional divisors of a crepant resolution of
$\tilde X/K$ equals the number of irreducible components of
the fixed point locus of $K$ on $\tilde\calX$.
}
NumbEx%
\finishproclaim
One can check that $K$ contains 6 elements which have nodes
as  isolated fixed points. Each of them fixes 16 nodes.
So each node occurs as fixed point of $K$. Hence all 96 exceptional
lines on $\tilde X$ are in the fixed point locus of $K$.
Hence we have to count only the one dimensional fixed curves
in $\calX$. This can be done withe results of the previous
section. We just give the result: There are 12 elements of $K$ having
a one dimensional fixed point locus and each of them has 4 components,
which are elliptic curves.
\proclaim
{Lemma}
{The number of components of the fixed point locus of $K$ on
$\tilde\calX$ is $144=96+12\cdot 4$.}
ResFix%
\finishproclaim
Now we are able ro compute the Picard number of a Calabi-Yau model
of $\calX/K$. The Picard number of the regular locus can be
computed by means of the results of section 6, especially
theorem \DecIrr. The result of a computation is 4. Hence we get:
\proclaim
{Lemma}
{The Picard number of a Calabi-Yau model of $X(\Gamma)$,
$\Gamma=\Gamma_{2,0}[4]\cap\Gamma_2[2]$, is $148$.}
PicLas%
\finishproclaim
We want to compute now the Euler number. Since $K$ is abelian, the string
theoretic formula gives
$$\eqalign{
e&={1\over 32}\sum_{(g,h)\in K\times K}e(\tilde\calX^{<g,h>})\cr
&={1\over 32}\cdot64+{3\over 32}\sum_{g\ne \id}e(\tilde\calX^{g})
+{1\over 32}\sum_{\id\ne g\ne h\ne\id}e(\tilde\calX^{<g,h>}).\cr
\cr}$$
Since the fixed point set of a single involution is
an elliptic curve or one of the 96 exceptional lines, we get
$$e=20+{1\over 32}\sum_{\id\ne g\ne h\ne\id}e(\tilde\calX^{<g,h>}).$$
\vskip1cm\noindent
we still have to discuss how for two different $g,h$,
which are different from the identity, the fixed point loci
intersect in $\tilde\calX$. We want to compare this with the intersection
of the fixed point loci on the singular model $\calX$.
We have to discuss two cases, where the fixed point loci of
$g$ and $h$ both are curves in $\calX$ or where the fixed point
locus of one of them consists of nodes.
We start with the case:
\smallni
{\it First case.\/} $g$ fixes a curve in $\calX$ and $h$
fixes a node.
\smallni
There are 12 $g$ which fix a curve and 6 $h$ which fix a node.
Hence we have 72 cases to consider. In 48 cases the intersection
of the fixed point loci in $\calX$ is empty.
Hence only 24 pairs are of interest. In each case
the fixed locus $\Fix(g)$ of $g$ is the union of 4
smooth elliptic curves
$$\Fix(g)=E_1\cup E_2\cup E_3\cup E_4.$$
and the fixed point locus of of $h$ consists of 16 nodes.
The intersection of $\Fix(g)$ and $\Fix(h)$ consists of 8 nodes.
Each single $E_i$ contains 4 of these 8 nodes.
This shows that in each of the 8 nodes two of the 4 elliptic
curves come together. Now we consider $\tilde\calX$.
Since the fixed point set of $g$ is smooth, it consists of four
elliptic curves $\tilde E_1,\dots,\tilde E_4$, such that
the natural projection $\tilde E_i\to E_i$ is biholomorphic.
Let $a$ be one of the 8 nodes in $\Fix(g)\cap\Fix(h)$.
We can assume that $E_1,E_2$ are the two elliptic curves which
run into $a$.
Let $C$ be the
exceptional line over $a$.
Then $g$ induces an automorphism
of $C$ of order two. Since an involution $P^1$ has two fixed
points, we see that $\tilde E_1$ and $\tilde E_2$
each hit $C$ in one intersection point and both points are different.
So each of the 8 exceptional lines carries two
intersection points. This shows:
\proclaim
{Lemma}
{Let $g\in K$ be an element with a one dimensional fixed point set,
and $h\in K$ an element, which fixes nodes. There are $24$
possibilities. The joint fixed point locus on $\tilde\calX$
consists of $16$ points.}
gFixh%
\finishproclaim
In the formula for the Euler number each pair $(g,h)$ of the above
form contributes with $16/32$. We have 24 pairs. Together with the
pairs $(h,g)$ we get the contribution $24$ to the Euler number.
Hence we have
$$e=44+ {1\over 32}\sum_{{\id\ne g\ne h\ne\id\atop
\dim\Fix(g)=\dim\Fix(h)=1}}e(\tilde\calX^{<g,h>})$$
{\it Second case.\/} Both $g$ and $h$ have one dimensional fixed
point locus on $\calX$ (4 elliptic curves).
The number of intersection points of $\Fix(g)$ and $\Fix(h)$ on
$\calX$ is 0, 8 or 16. The number of pairs $(g,h)$ with 8 intersection points
is 24 and that with 16 intersection points is 48.
\smallni
{\it Pairs with $16$ intersection points.\/}
In this case one can check that none of the 16 is a node, and one can
check furthermore that Hence the contribution to the Euler
for each such pair is $(1/32)\cdot 16=1/2$.
\smallni
{\it Pairs with $8$ intersection points.\/}
In this case one can check that all 8 intersection points are nodes.
Let $a$ be such a node. One can see that that two of the components of
$\Fix(g)$ run into $a$ and the same is true for $\Fix(h)$.
Hence as in the first case above $g$ has two fixed points $a_1,a_2$
on the exceptional line $C$ over $a$ and $h$ has also two fixed
points $b_1,b_2$ there. Now $g$ and $h$ induce involutions, which
commute. But commuting involutions on $P^1$ don't have a joint fixed point, since they
alway can be transformed to the involutions $z\mapsto-z$ and $z\mapsto-1/z$
on the Riemann sphere.
Hence there is no contribution to the Euler number.
We get as  contribution
$24$ to the Euler number. This gives
$$e=68$$
for the Euler number.
\proclaim
{Theorem}
{A Calabi-Yau model which belongs to $\Gamma_{2,0}[4]\cap\Gamma_2[2]$
has Hodge numbers $h^{11}=148$, $h^{12}= 114$.}
GamDia%
\finishproclaim
\bye